\documentstyle[11pt,pstricks,pst-node,amssymb]{article}

\newcommand{\qed}{\relax{\ifhmode\unskip\nobreak\hfill$\Box$\fi
	\ifmmode\ifinner\else\hskip5pt\fi \hfill\Box\fi}\relax}
\newcommand{\ie}{{\em i.e.}}

\newcommand{\mto}{\lput{:U}{\pspicture(0,0)(0,0)
\psline[arrows=->,arrowscale=1.5](2.2pt,0)(2.3pt,0)\endpspicture}}
\newcommand{\mfro}{\lput{:U}{\pspicture(0,0)(0,0)
\psline[arrows=->,arrowscale=1.5](-2.2pt,0)(-2.3pt,0)\endpspicture}}

\newcommand{\sixvertex}[6]{
\pspicture[.5](-.7,-1.7)(.7,.7)
\pcline(0,0)(0,.7)\mto    \pcline(0,0)(.7,0)\mfro
\pcline(0,0)(-.7,0)\mfro  \pcline(0,0)(0,-.7)\mto
\rput(0,-1.2){$#1$}\endpspicture\hspace{1cm}
\pspicture[.5](-.7,-1.7)(.7,.7)
\pcline(0,0)(0,.7)\mfro   \pcline(0,0)(.7,0)\mto
\pcline(0,0)(-.7,0)\mto   \pcline(0,0)(0,-.7)\mfro
\rput(0,-1.2){$#2$}\endpspicture\hspace{1cm}
\pspicture[.5](-.7,-1.7)(.7,.7)
\pcline(0,0)(0,.7)\mto    \pcline(0,0)(.7,0)\mto
\pcline(0,0)(-.7,0)\mfro  \pcline(0,0)(0,-.7)\mfro
\rput(0,-1.2){$#3$}\endpspicture\hspace{1cm}
\pspicture[.5](-.7,-1.7)(.7,.7)
\pcline(0,0)(0,.7)\mfro   \pcline(0,0)(.7,0)\mfro
\pcline(0,0)(-.7,0)\mto   \pcline(0,0)(0,-.7)\mto
\rput(0,-1.2){$#4$}\endpspicture\hspace{1cm}
\pspicture[.5](-.7,-1.7)(.7,.7)
\pcline(0,0)(0,.7)\mto    \pcline(0,0)(.7,0)\mfro
\pcline(0,0)(-.7,0)\mto   \pcline(0,0)(0,-.7)\mfro
\rput(0,-1.2){$#5$}\endpspicture\hspace{1cm}
\pspicture[.5](-.7,-1.7)(.7,.7)
\pcline(0,0)(0,.7)\mfro   \pcline(0,0)(.7,0)\mto
\pcline(0,0)(-.7,0)\mfro  \pcline(0,0)(0,-.7)\mto
\rput(0,-1.2){$#6$}\endpspicture}
\newcounter{fignum}
\newcommand{\df}{\em}
\def\bqn {\begin{equation}}
\def\eqn {\end{equation}}

\newtheorem{theorem}{Theorem}
\newtheorem{lemma}[theorem]{Lemma}
\newenvironment{proof}{{\sc Proof: }}{\qed}

\psset{linewidth=.25pt,dash=3pt 3pt}
\newrgbcolor{lightblue}{.5 .5 1}
\SpecialCoor

\begin{document}

\title{Another proof of the alternating sign matrix conjecture}
\author{Greg Kuperberg}
\date{December 10, 1995}

\maketitle

\begin{abstract}
Mills, Robbins, and Rumsey \cite{MRR:asm} conjectured, and Zeilberger
\cite{Zeilberger:asm} recently proved, that there are ${1!4!7!\ldots
(3n-2)!\over n!(n+1)!\ldots (2n-1)!}$ alternating sign matrices of
order $n$.  We give a new proof of this result using an analysis of the
six-vertex state model (also called square ice) based on the
Yang-Baxter equation.
\end{abstract}

Mills, Robbins, and Rumsey~\cite{MRR:asm} conjectured that:

\begin{theorem}[Zeilberger] There are 
$$A(n) = {1!4!7!\ldots (3n-2)! \over n!(n+1)!(n+2)!\ldots (2n-1)!}$$
$n \times n$ alternating sign matrices. \label{thmain}
\end{theorem}

Here, an {\df alternating sign matrix} or {\df ASM} is a matrix of
$0$'s, $1$'s, and $-1$'s such that the non-zero elements in each row
and column alternate between $1$ and $-1$ and begin and end with $1$,
for example:
$$\left(\begin{array}{rrrr}
 0 & 1 & 0 & 0 \\
 1 &-1 & 1 & 0 \\
 0 & 0 & 0 & 1 \\
 0 & 1 & 0 & 0
\end{array}\right)$$
Alternating sign matrices are related to a number of other
combinatorial objects that, remarkably, are also enumerated or
conjectured to be enumerated by ratios of progressions of factorials or
staggered factorials~\cite{Robbins:intelligencer,Stanley:symmetries}.

Zeilberger~\cite{Zeilberger:asm} recently proved Theorem~\ref{thmain} by
establishing that ASM's are equinumerous with totally symmetric,
self-complementary plane partitions, which were enumerated by
Andrews~\cite{Andrews:tsscpp}.  In this paper, we present a new proof.  The
most interesting part of the proof is due to Izergin and
Korepin~\cite{Izergin:six-vertex,Korepin:qism}, who follow Baxter's
remarkable use of the Yang-Baxter equation~\cite{Baxter:exactly}.

If $x$ is a number, define the $x$-enumeration $A(n;x)$ of $n \times n$ ASM's
as their total weight, where the weight of an individual matrix is $x^k$ if
it has $k$ entries equal to $-1$. A variation of the proof
establishes another conjecture of Mills, Robbins, and Rumsey:

\begin{theorem} ASM's are 3-enumerated by
$$A(2n+1;3) = \left(3^{n(n+1)/2}{2!5!8!\ldots(3n-1)! \over
(n+1)!(n+2)!\ldots(2n)!}\right)^2$$
$$A(2n;3) = 3^{n-1} {(3n-1)!(n-1)! \over (2n-1)!^2}A(2n-1;3)$$
\label{thmain3}
\end{theorem}

A second variation establishes the well-known 2-enumeration
of ASM's \cite{EKLP,MRR:asm}:
$$A(n;2) = 2^{n(n-1)/2}.$$
Finally, the following result, also conjectured by Mills, Robbins,
and Rumsey, follows easily from the general method:

\begin{theorem} For each $n$, there exists a polynomial $B(n;x)$
such that
$$A(n;x) = B(n;x)B(n+1;x)$$
for $n$ odd and
$$A(n;x) = 2B(n;x)B(n+1;x)$$
for $n$ even. \label{thx}
\end{theorem}

Mills, Robbins, and Rumsey further conjectured that for $n$ odd, $B(n;x)$ is
the $x$-enumeration of vertically symmetric ASM's (where the weight is $x^k$
if there are $k$ ones to the left of the middle column), but this relation
remains open.

\subsection{Acknowledgements}

The author would like to thank David Robbins and Doron Zeilberger for
explaining the significance of the alternating sign matrix conjecture, Igor
Frenkel and Vaughan Jones for explaining the Yang-Baxter equation, and Jim
Propp for introducing the author to enumerative combinatorics.  Also, the
computer algebra package Maple~\cite{Maple} was indispensible for
computations, and the \TeX\, macro package PSTricks~\cite{PSTricks} was
very useful for setting the figures.

\section{State sums}

The {\df six-vertex model} in general
refers to the multiplicative weighted enumeration of orientations
of a tetravalent planar graph $G$ (called {\df states})
such that at each vertex, two arrows go in and two go out.  Number
the six allowed orientations incident to a given vertex (called
states of a vertex) 1 through 6:
$$
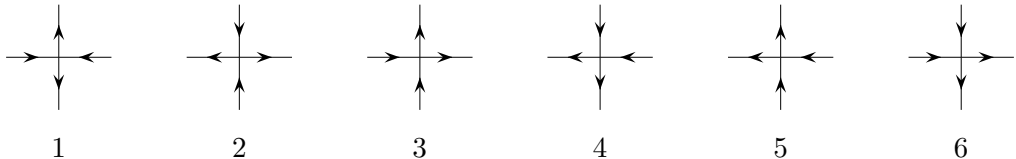
\sixvertex{1}{2}{3}{4}{5}{6}
$$
State $i$ at vertex $v$ is given a weight $w(i,v)$. The weight of a state of
$G$ is the product of the weights of its vertices, and the {\df state sum} is
the total weight of all states. The six-vertex model may also be considered
with boundary conditions, meaning that there may be univalent vertices whose
edges have fixed orientations.  In particular, consider a six-vertex state of
an $n \times n$ square grid with edges pointing inward at the sides and
outward at the top and bottom:
$$
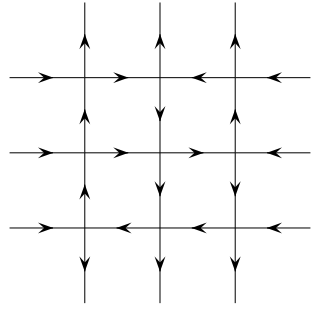
\pspicture(-2,-2)(2,2)
\pnode(-2, 2){aa}\pnode(-1, 2){ab}\pnode( 0, 2){ac}\pnode(1, 2){ad}\pnode(2, 2){ae}
\pnode(-2, 1){ba}\pnode(-1, 1){bb}\pnode( 0, 1){bc}\pnode(1, 1){bd}\pnode(2, 1){be}
\pnode(-2, 0){ca}\pnode(-1, 0){cb}\pnode( 0, 0){cc}\pnode(1, 0){cd}\pnode(2, 0){ce}
\pnode(-2,-1){da}\pnode(-1,-1){db}\pnode( 0,-1){dc}\pnode(1,-1){dd}\pnode(2,-1){de}
\pnode(-2,-2){ea}\pnode(-1,-2){eb}\pnode( 0,-2){ec}\pnode(1,-2){ed}\pnode(2,-2){ee}
\ncline{ba}{bb}\mto \ncline{bb}{bc}\mto \ncline{bc}{bd}\mfro\ncline{bd}{be}\mfro
\ncline{ca}{cb}\mto \ncline{cb}{cc}\mto \ncline{cc}{cd}\mto \ncline{cd}{ce}\mfro
\ncline{da}{db}\mto \ncline{db}{dc}\mfro\ncline{dc}{dd}\mfro\ncline{dd}{de}\mfro
\ncline{ab}{bb}\mfro\ncline{bb}{cb}\mfro\ncline{cb}{db}\mfro\ncline{db}{eb}\mto
\ncline{ac}{bc}\mfro\ncline{bc}{cc}\mto \ncline{cc}{dc}\mto \ncline{dc}{ec}\mto
\ncline{ad}{bd}\mfro\ncline{bd}{cd}\mfro\ncline{cd}{dd}\mto \ncline{dd}{ed}\mto
\endpspicture
$$
The six-vertex model on a square grid is also called {\df square ice}. A
square ice state can be converted to an ASM by the
correspondence
$$
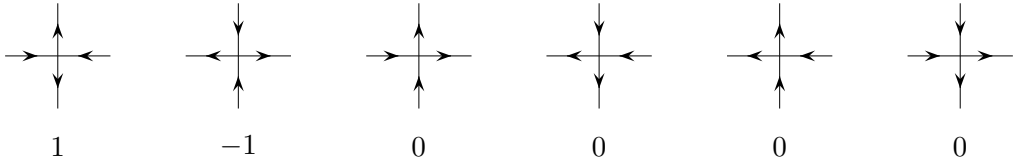
\sixvertex{1}{-1}{0}{0}{0}{0}
$$
This conversion is bijective \cite{EKLP,RR:asm}.  Thus, the enumeration of
ASM's is equivalent to a six-vertex state sum in which all weights are 1.

Let $h$ be a complex number or an indeterminate, let $q^x$ denote $e^{hx}$,
and let $[x]$ denote ${q^{x/2} - q^{-x/2} \over q^{1/2} - q^{-1/2}}$.  We
will consider various {\df half-integral Laurent polynomials}, meaning
polynomials with integral or half-integral exponents of either sign such that
the difference between any two exponents is an integer.  For example, if $q$
is fixed, $[x]$ is a half-integral Laurent polynomial in $q^x$.  Given two
such polynomials $P(t)$ and $Q(t)$ over a ring $A$, we will say that $Q$
divides $P$ if $P(t)/Q(t) \in A[t^{1/2},t^{-1/2}]$. For example, $t$
divides $1$.

A vertex labelled by $x$:
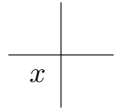
$$
\pspicture(-.7,-.7)(.7,.7)
\psline(0,0)(0,.7)    \psline(0,0)(.7,0)
\psline(-0,0)(-.7,0)  \psline(0,-0)(0,-.7)
\rput[tr](-.2,-.2){$x$}\endpspicture
$$
denotes the six weights:
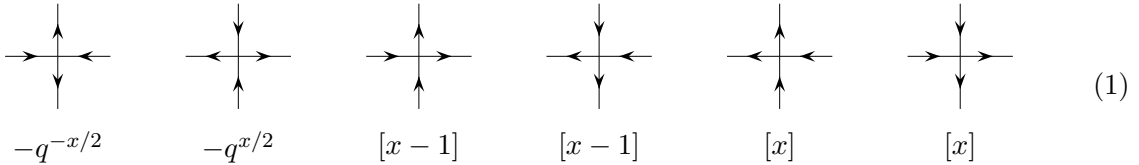
\bqn
\sixvertex{-q^{-x/2}}{-q^{x/2}}{[x-1]}{[x-1]}{[x]}{[x]} \label{fweights}
\eqn
(Since the weights are invariant under rotation by 180 degres,
but not 90 degrees, the meaning of a vertex depends on which
pair of kitty-corner quadrants contains its label.)
Such a vertex is called an $R$-matrix and is also denoted as $R(x)$. 

\begin{theorem}[Baxter]  If $x
= y + z$,  the $R$-matrices $R(x)$, $R(y)$, and $R(z)$ satisfy the
equation
$$
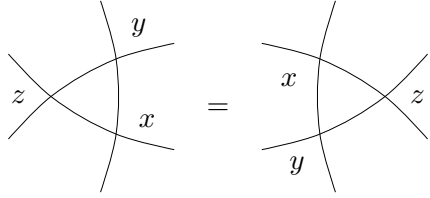
\pspicture[.42](-1.8,-1)(1.1,2)
\pnode(0,0){b1}\pnode(0,1){b2}\pnode(-.866,.5){b3}
\pnode([angle=255,nodesep=.8]b1){a1}\pnode([angle=345,nodesep=.8]b1){a2}
\pnode([angle= 15,nodesep=.8]b2){a3}\pnode([angle=105,nodesep=.8]b2){a4}
\pnode([angle=135,nodesep=.8]b3){a5}\pnode([angle=225,nodesep=.8]b3){a6}
\pscurve(a1)(b1)(b2)(a4)
\pscurve(a3)(b2)(b3)(a6)
\pscurve(a5)(b3)(b1)(a2)
\rput[bl](.3,.1){$x$}   % x = y + z
\rput[bl](.2,1.3){$y$}
\rput[r](-1.2,.5){$z$}
\endpspicture = 
\pspicture[.42](-1.1,-1)(1.8,2)
\pnode(0,0){b1}\pnode(0,1){b2}\pnode(.866,.5){b3}
\pnode([angle=285,nodesep=.8]b1){a1}\pnode([angle=195,nodesep=.8]b1){a2}
\pnode([angle=165,nodesep=.8]b2){a3}\pnode([angle= 75,nodesep=.8]b2){a4}
\pnode([angle= 45,nodesep=.8]b3){a5}\pnode([angle=315,nodesep=.8]b3){a6}
\pscurve(a1)(b1)(b2)(a4)
\pscurve(a3)(b2)(b3)(a6)
\pscurve(a5)(b3)(b1)(a2)
\rput[tr](-.2,-.3){$y$}
\rput[tr](-.3,.8){$x$}
\rput[l](1.2,.5){$z$}
\endpspicture
$$
\end{theorem}
This remarkable identity is known as the {\df star-triangle} relation or the
{\df Yang-Baxter} equation\cite{Baxter:exactly}. Specifically, $R(x)$ is said
to parameterize the {\df trigonometric solutions} to the Yang-Baxter
equation. Before proving it, we discuss exactly what the equation means. Each
of the two graphs in the equation has six external edges, meaning edges with
a univalent vertex.  For each external edge on the left, there is a
corresponding external edge on the right whose univalent vertex is in the
same position; for example, on both sides there is a lowest univalent
endpoint, and the two edges with this endpoint correspond to each other.  For
each of the 64 orientations of the external edges on the left, one can form a
state sum $Z$ by summing over admissible orientations of the three internal
edges, and one can consider the same orientation on the right and form
another state sum $Z'$. The equation then says that the $Z=Z'$ in all 64
cases. In order for the state sum to be non-zero, three edges must point in
and three must point out, so the identity is trivial in 44 of the 64 cases. 
Note further that the equation simply says that the left side is invariant
under rotation by 180 degrees, so the 20 non-trivial numerical identities
reduce to 10 identities repeated twice.  The argument that follows uses other
tricks to further reduce the number of numerical identities to one which can
checked be checked easily:

\begin{proof} We first rearrange the left side of the Yang-Baxter equation:
$$
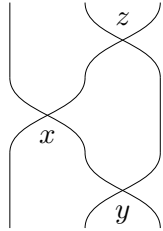
\pspicture(0,0)(2,3)
\qline(0,0)(0,1)
\pccurve[angleA=90,angleB=270](1,0)(2,1)
\pccurve[angleA=90,angleB=270](2,0)(1,1)
\qline(2,1)(2,2)
\pccurve[angleA=90,angleB=270](0,1)(1,2)
\pccurve[angleA=90,angleB=270](1,1)(0,2)
\qline(0,2)(0,3)
\pccurve[angleA=90,angleB=270](1,2)(2,3)
\pccurve[angleA=90,angleB=270](2,2)(1,3)
\rput[t](1.5,.3){$y$}
\rput[t](.5,1.3){$x$}
\rput[b](1.5,2.7){$z$}
\endpspicture
$$
Consider the
following augmentation of the six-vertex model:  Suppose that a graph has a
curved edge with a horizontal tangent at a point $p$ and which is
concave down at $p$.  If the edge is oriented to the left in some six-vertex
state, $p$ is assigned a multiplicative weight of $-q^{1/2}$, but if it
points to the right, it is assigned a multiplicative weight of 1.
Contrariwise, if the tangent is horizontal but the curve is concave
up, $p$ has weight $-q^{-1/2}$ when the edge points to the left
and weight $1$ when it points to the right.  With this convention,
the following simple identities hold:
\begin{equation}
\pspicture[.4](-1,-1)(1,1)
\qline(-.7,-1)(-.7,0)
\psarc(-.35,0){.35}{0}{180}\psarc(.35,0){.35}{180}{0}
\qline(.7,0)(.7,1)
\endpspicture
=
\pspicture[.4](-1.2,-1)(1.2,1)
\qline(.7,-1)(.7,0)
\psarc(.35,0){.35}{0}{180}\psarc(-.35,0){.35}{180}{0}
\qline(-.7,0)(-.7,1)
\endpspicture
=
\pspicture[.4](-.5,-.7)(.5,.7)
\qline(0,-.7)(0,.7)
\endpspicture
\hspace{2cm}
\pspicture[.4](-.7,-.5)(.7,.5)
\pscircle(0,0){.5}
\endpspicture
= -q^{1/2} - q^{-1/2} = -[2]
\label{etl}
\end{equation}
Moreover, $R(x)$ can be expressed as
$$
\pspicture[.4](-.7,-.7)(.7,.7)
\psline(.7;45)(.7;225)\psline(.7;135)(.7;315)
\rput[t](0,-.2828){$x$}\endpspicture
= 
[x]
\pspicture[.4](-.9,-.7)(.9,.7)
\pccurve[angleA= 45,angleB=135,ncurv=1,nodesep=0](.7;225)(.7;315)
\pccurve[angleA=225,angleB=315,ncurv=1,nodesep=0](.7; 45)(.7;135)
\endpspicture
+[x-1]
\pspicture[.4](-.9,-.7)(.9,.7)
\pccurve[angleA=225,angleB=135,ncurv=1,nodesep=0](.7; 45)(.7;315)
\pccurve[angleA= 45,angleB=315,ncurv=1,nodesep=0](.7;225)(.7;135)
\endpspicture
$$
Thus, a six-vertex state sum involving $R$-matrices can be expanded as a sum
of curves in a calculus in which each closed loop contributes a factor of
$-[2]$. (This calculus is called the {\df Temperley-Lieb category} and is
closely related to the quantum group $U_q(\mbox{sl}(2))$
\cite{KL:recoupling}.) The calculus is invariant under isotopy of curves by
equation~(\ref{etl}). The left side of the Yang-Baxter equation then expands
to eight terms, which may collected into five terms corresponding to the five
crossingless matchings of six points on a circle. Three of the matchings are
invariant under rotation by 180 degrees.  The coefficients of the other two are
$$[z-1][x][y-1]$$
and
$$[z][x-1][y-1] + [z][x][y] + [z-1][x-1][y] - [2][z][x-1][y].$$
These two quantities are rendered equal by the identities
$x = y+z$, $[-a] = -[a]$, and 
$$[a][b] - [a+1][b-1] = [a-b+1].$$
Thus, the left side is invariant under rotation by 180 degrees.
\end{proof}

As a final notational convenience, define
$$
\pspicture[.5](-.9,-1.2)(.9,.9)
\psline(0,0)(0,.7)    \psline(0,0)(.7,0)
\psline(-0,0)(-.7,0)  \psline(0,-0)(0,-.7)
\rput[r](-.9,0){$x$}
\rput[t](0,-.9){$y$}
\endpspicture
=
\pspicture[.5](-1.6,-1.2)(.7,.9)
\psline(0,0)(0,.7)    \psline(0,0)(.7,0)
\psline(-0,0)(-.7,0)  \psline(0,-0)(0,-.7)
\rput[rt](-.2,-.2){$x-y$}
\endpspicture
$$
when the lines rather than the vertices of a tetravalent graph are labelled.
Following Izergin and Korepin~\cite{Izergin:six-vertex,Korepin:qism},
consider $n \times n$ square ice with arbitrary parameters $X =
x_0,\ldots,x_{n-1}$ and $Y = y_0,\ldots,y_{n-1}$ for the horizontal and
vertical lines:
$$
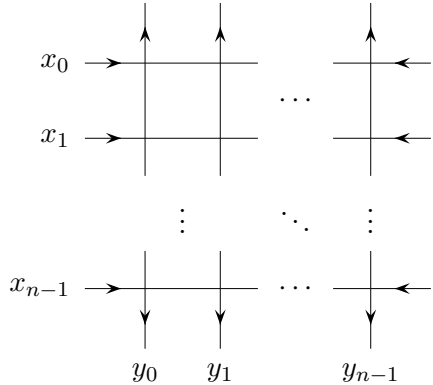
\pspicture(-1.5,-1.5)(4,4)
\pnode(-1, 4){aa}\pnode(0, 4){ab}\pnode(1, 4){ac}\pnode(1.5, 4){ad}
\pnode(2.5, 4){ax}\pnode(3, 4){ay}\pnode(4, 4){az}
\pnode(-1, 3){ba}\pnode(0, 3){bb}\pnode(1, 3){bc}\pnode(1.5, 3){bd}
\pnode(2.5, 3){bx}\pnode(3, 3){by}\pnode(4, 3){bz}
\pnode(-1, 2){ca}\pnode(0, 2){cb}\pnode(1, 2){cc}\pnode(1.5, 2){cd}
\pnode(2.5, 2){cx}\pnode(3, 2){cy}\pnode(4, 2){cz}
\pnode(-1,1.5){da}\pnode(0,1.5){db}\pnode(1,1.5){dc}\pnode(1.5,1.5){dd}
\pnode(2.5,1.5){dx}\pnode(3,1.5){dy}\pnode(4,1.5){dz}
\pnode(-1,.5){xa}\pnode(0,.5){xb}\pnode(1,.5){xc}\pnode(1.5,.5){xd}
\pnode(2.5,.5){xx}\pnode(3,.5){xy}\pnode(4,.5){xz}
\pnode(-1, 0){ya}\pnode(0, 0){yb}\pnode(1, 0){yc}\pnode(1.5, 0){yd}
\pnode(2.5, 0){yx}\pnode(3, 0){yy}\pnode(4, 0){yz}
\pnode(-1,-1){za}\pnode(0,-1){zb}\pnode(1,-1){zc}\pnode(1.5,-1){zd}
\pnode(2.5,-1){zx}\pnode(3,-1){zy}\pnode(4,-1){zz}
\ncline[nodesepA=.2]{ba}{bb}\mto \ncline{bb}{bc}\ncline{bc}{bd}
\ncline{bx}{by}\ncline[nodesepB=.2]{by}{bz}\mfro
\ncline[nodesepA=.2]{ca}{cb}\mto \ncline{cb}{cc}\ncline{cc}{cd}
\ncline{cx}{cy}\ncline[nodesepB=.2]{cy}{cz}\mfro
\ncline[nodesepA=.2]{ya}{yb}\mto \ncline{yb}{yc}\ncline{yc}{yd}
\ncline{yx}{yy}\ncline[nodesepB=.2]{yy}{yz}\mfro
\ncline[nodesepA=.2]{ab}{bb}\mfro\ncline{bb}{cb}\ncline{cb}{db}
\ncline{xb}{yb}\ncline[nodesepB=.2]{yb}{zb}\mto 
\ncline[nodesepA=.2]{ac}{bc}\mfro\ncline{bc}{cc}\ncline{cc}{dc}
\ncline{xc}{yc}\ncline[nodesepB=.2]{yc}{zc}\mto 
\ncline[nodesepA=.2]{ay}{by}\mfro\ncline{by}{cy}\ncline{cy}{dy}
\ncline{xy}{yy}\ncline[nodesepB=.2]{yy}{zy}\mto 
\rput[r](ba){$x_0$}\rput[r](ca){$x_1$}\rput[r](ya){$x_{n-1}$}
\rput[t](zb){$y_0$}\rput[t](zc){$y_1$}\rput[t](zy){$y_{n-1}$}
\rput(.5,1){$\vdots$}\rput(2,1){$\ddots$}\rput(2,2.5){$\cdots$}
\rput(2,0){$\cdots$}\rput(3,1){$\vdots$}
\endpspicture
$$
Let $Z(n;X,Y)$ be the resulting state sum.

\begin{lemma}[Baxter] The function $Z(n;X,Y)$ is symmetric in the $x_i$'s
and in the $y_i$'s.
\label{lsymmetric}
\end{lemma}
\begin{proof} Consider the $i$th and $i+1$st horizontal lines. An extra vertex
(implicitly labelled by $x_i-x_{i+1}$) may be introduced on the left at the
expense of a generically non-zero multiplicative factor:
$$\pspicture[.42](-1.5,-1)(2.5,2)
\pcline(-.8,0)(-0,0)\mto\pcline(-.8,1)(-0,1)\mto
\psline(0,0)(1.5,0)\psline(0,1)(1.5,1)
\psline(0,-.5)(0,1.5)\psline(1,-.5)(1,1.5)
\rput(2,.5){$\cdots$}
\rput[r](-1,0){$x_{i+1}$}
\rput[r](-1,1){$x_i$}
\endpspicture
= [x_i-x_{i+1}-1]
\pspicture[.42](-2.8,-1)(2.5,2)
\pnode(0,0){a1}\pnode(0,1){a2}\pnode(-.866,.5){b1}
\pnode(-1.732,0){a3}\pnode(-1.732,1){a4}
\nccurve[angleA=180,angleB=-30]{a1}{b1}\nccurve[angleA=180,angleB=30]{a2}{b1}
\ncline{a3}{b1}\mto\ncline{a4}{b1}\mto
\psline(0,0)(1.5,0)\psline(0,1)(1.5,1)
\psline(0,-.5)(0,1.5)\psline(1,-.5)(1,1.5)
\rput(2,.5){$\cdots$}
\rput[r](-1.932,0){$x_i$}
\rput[r](-1.932,1){$x_{i+1}$}
\endpspicture
$$
This relation holds because in an allowed state, all four edges of the new
vertex must point to the right. By the Yang-Baxter equation, the vertex can
be moved from the left side to the right, whereupon it can be removed, which
recovers the multiplicative factor.  This operation switches the labels $x_i$
and $x_{i+1}$. Therefore $Z(n;X,Y)$ is symmetric in $x_i$ and $x_{i+1}$ for
each $i$, which renders it symmetric in all $x_i$'s.  The same argument
applies to the $y_i$'s.
\end{proof}

\begin{lemma} If $x_i = y_j + 1$, then 
$$Z(n;X,Y) = -q^{-1/2}\left(\prod_{k \ne i} [x_i-y_k]\right)
\left(\prod_{k \ne j} [x_k - y_j]\right)
Z(n-1;X \setminus x_i, Y \setminus y_j).$$
\label{lzrelation}
\end{lemma}
\begin{proof}
Assume first that $i=j=0$.  By Figure~(\ref{fweights}), the upper left vertex
must have state 1 in a non-zero state of the grid. This forces the rest of
the top row to have state 5 and the rest of the left column to have state 6,
which yields the given multiplicative factor. (In terms of ASM's, only those
matrices with a $1$ in the top left corner contribute.)  The remainder of the
grid is an $n-1 \times n-1$ square ice state.

The general case follows from Lemma~\ref{lsymmetric}.
\end{proof}

\begin{lemma} The quantity $q^{nx_0/2}Z(n;X,Y)$ is a polynomial in
$q^{x_0}$ of degree at most $n-1$. \label{lzdegree}
\end{lemma}
\begin{proof}
If we multiply all weights of vertices in the first row by $q^{x_0/2}$,
then $q^{x_0}$ appears linearly in those weights in which it appears 
at all.  Therefore the modified state sum
$$Z'(n) = q^{nx_0/2}Z(n;X,Y)$$
is a polynomial in $q^{x_0}$.  The first row is the only row in which $x_0$
appears.  In this row, there must be one vertex in state 1, whose modified
weight does not involve $x_0$, and $n-1$ vertices in state 5 or 3.  (In terms
of ASM's, there must be a 1 in the top row.) Therefore $Z'(n)$ has degree at
most $n-1$.
\end{proof}

\begin{theorem}[Izergin,Korepin] The state sum $Z(n;X,Y)$ is given by
$$Z(n;X,Y) = {(-1)^n\left(\prod_{i=0}^{n-1} q^{(y_i-x_i)/2}\right)
\prod_{0 \le i,j < n} [x_i-y_j][x_i-y_j-1] \over \left(\prod_{0 \le j<i<n}
[x_i-x_j]\right)\left(\prod_{0 \le i < j < n} [y_i-y_j]\right)}\det M,$$
where
$$M_{i,j} = {1 \over [x_i - y_j][x_i - y_j - 1]}.$$
\label{thkorepin}
\end{theorem}
\begin{proof}
Lemmas~\ref{lzrelation} and \ref{lzdegree}, together with $Z(0) = 1$,
inductively determine $Z(n)$ by Lagrange interpolation.  It is routine to
check that the right side satisfies Lemma~\ref{lzrelation}.  To check that it
also satisfies Lemma~\ref{lzdegree}, Let $P$ be the numerator, let $Q$ be the
denominator, let $D$ be the determinant, and let $D'$ be a term in the
expansion of the determinant.  The product $PD$ is a half-integral Laurent
polynomial because $PD'$ is for any choice of $D'$.  Moreover, $Q$ divides
$PD$, because $D$ is antisymmetric in the $x_i$'s and in the $y_j$'s and
therefore in the $q^{x_i}$'s and in the $q^{y_j}$'s. Thus, ${P\over Q}D$ is a
half-integral Laurent polynomial polynomial in $q^{x_0}$.  Finally, the
leading term (expanded as a Laurent polynomial in $q^{x_0}$) of any $PD'$ has
exponent $(2n-3)/2$, while the trailing term has exponent $(1-2n)/2$. 
Therefore the same is true of $PD$, and ${P \over Q}D$ has leading exponent
at most $(n-2)/2$ and trailing exponent at least $-n/2$.  In conclusion,
$q^{nx_0/2} {P \over Q}D$ is a polynomial in $q^{x_0}$ and has degree at
most $n-1$.
\end{proof}

\section{Determinants}

Consider the state-sum value
$$Z_{\frac12}(n) = Z(n;\frac12,\frac12,\ldots,\frac12,0,0,\ldots,0).$$
In any $n \times n$ square ice state, there are $n$ more vertices
in state 1 than state 2, equal numbers in states 3 and 4, and equal
numbers in states 5 and 6.  Since the weights of these states
in $R(2)$ are $-q^{-1/4}$, $-q^{1/4}$, $-[\frac12]$, $-[\frac12]$, $[\frac12]$,
and $[\frac12]$, respectively, it follows that
\bqn
A(n;x) = [\frac12]^{n-n^2}(-1)^n q^{n/4} Z_{\frac12}(n),\label{ean}
\eqn
where $x = 1/[\frac12]^2 = [2]+2$.  Unfortunately, the determinant
in Theorem~\ref{thkorepin} is singular for $Z_{1/2}(n)$.  Therefore,
we will instead evaluate
$$Z_{\frac12}(n;\epsilon) =
Z(n;\frac12+\epsilon,\frac12+2\epsilon,\ldots,\frac12+n\epsilon,0,
-\epsilon,-2\epsilon\ldots,(1-n)\epsilon)$$
when $h = {4\pi \sqrt{-1} \over 3}$, which implies that $x=1$.

Let $s = q^{\epsilon}$.  Firstly,
$$[k\epsilon+\frac12][k\epsilon-\frac12] = {s^k + 1 + s^{-k} \over -3}$$
and
$$[k\epsilon] = {s^{k/2} - s^{-k/2} \over \sqrt{-3}}.$$
The matrix $M$ of Theorem~\ref{thkorepin}
becomes
$$M_{i,j} = {-3 \over s^{i+j+1}+ 1 + s^{-(i+j+1)}}.$$
The state sum becomes
$$Z_\frac12(n;\epsilon) = {q^{-n/4} s^{-n^2/2} 3^{-n(n+1)/2}
\prod_{0 \le i,j < n} (s^{i+j+1}+1+s^{-(i+j+1)})
\over \prod_{0 \le j < i < n} (s^{(i-j)/2}-s^{(j-i)/2})^2} \det M.$$
The determinant of $M$ can be computed using the following two lemmas, for
which we extend the bracket notation by defining $[x]_t = {t^{x/2} - t^{-x/2}
\over t - t^{-1}}$ for any $t$.

\begin{lemma}[Cauchy] Let $X = x_0,\ldots,x_{n-1}$ and $Y = y_0,\ldots,y_{n-1}$
be variables, and let
$$T(n;X,Y)_{i,j} = {1 \over [x_i-y_j]_t}$$
for $0 \le i,j < n$.  Then
$$\det T(n,k;t) = {\left(\prod_{0 \le j<i<n}
[x_i-x_j]_t\right)\left(\prod_{0 \le i < j < n} [y_i-y_j]_t\right)
\over \prod_{0 \le i,j < n} [x_i-y_j]_t}.$$
\label{ltdet}
\end{lemma}
\begin{proof}
Let $D$ be the determinant, let $P$ be the denominator, and let $Q$ be the
numerator.  Then the arguments of Theorem~\ref{thkorepin} apply, but with the
conclusion that ${P \over Q}D$ is a degree 0 polynomial in all variables,
\ie, a constant.  Let $D'$ be the diagonal term in the determinant; $D'$ is
the only term such that $PD'$ is not divisible by any $[x_i - y_i]_t$. All
factors of ${P \over Q}D'$ cancel at the specialization $x_i = y_i$ and all
other terms of ${P \over Q}D$ vanish; therefore ${P \over Q}D = 1$.
\end{proof}

Lemma~\ref{ltdet} can also be proved by induction
using Dodgson's condensation method~\cite{Robbins:intelligencer,RR:asm}.

Let $T(n) = T(n;1,2,\ldots,n,0,-1,-2,\ldots,1-n)$.  Then
$$T(n)_{i,j} = {1 \over [i+j+1]_t}$$
and $\det T(n)$ is given by Lemma~\ref{ltdet}.

\begin{lemma} Let
$$S(n;s,t)_{i,j} = {s^{(i+j+1)/2} - s^{-(i+j+1)/2} \over
t^{(i+j+1)/2} - t^{-(i+j+1)/2}}$$
for $0 \le i,j < n$.  Then
$$\det S(n;s,t) = {(-1)^{n(n-1)/2}\over (t^{1/2} - t{-1/2})^n}
\left(\prod_{0 \le j < i < n} [i-j]_t^2\right)
\left(\prod_{0 \le i,j < n} {s^{1/2}t^{(i-j)/2} - s^{1/2}t^{(j-i)/2}
\over [i+j+1]_t}\right).$$
\label{lsdet}
\end{lemma}
\begin{proof}
The quantity $s^{n^2/2}(\det S(n;s,t)$ has degree $n^2$ as a polynomial
in $s$.  Moreover, for $0 \le k < n$,
$$S(n;t^k,t)_{i,j} = [k]_{t^{i+j+1}} = \sum_{\ell = 0}^{k-1}
A(t^{\ell-(k-1)/2})_{i,j},$$
where the matrix $A(z)$ given by
$$A(z)_{i,j} = z^{i+j+1}$$
has rank 1.  Thus, the rank of $S(n;t^k,t)$ is at most $k$, and it follows
that $(s-t^k)^{n-k}$ divides $\det S(n;s,t)$. Similarly, $(s-t^{-k})^{n-k}$
divides the determinant for $1 \le k < n$.  These divisibilities determine
$\det S(n;s,t)$ up to a factor which is a function of $t$.  The leading
coefficient is then $(\det T(n))/(t^{1/2}-t^{-1/2})^n$.
\end{proof}

The determinant of $-M/3 = S(n;s,s^3)$ is given by Lemma~\ref{lsdet}.
Collecting factors yields
\begin{eqnarray*}
Z_{1/2}(n;\epsilon) & = &
{q^{-n/4} s^{-n^2/2} 3^{-n(n+1)/2} \left(\prod_{0 \le i,j < n}[3(i+j+1)]_s\right)
\over
\left(\prod_{0 \le i,j < n} [i+j+1]_s\right)
\left(\prod_{0 \le j < i < n} [i-j]^2_s\right)}
\\ & & {(-3)^n (-1)^{n(n-1)/2} \over [3]^n_s}
\left(\prod_{0 \le j < i < n} {[3(i-j)]^2_s\over [3]^2_s} \right)
\left(\prod_{0 \le i,j < n} {[3]_s[3(i-j)+1]_s \over [3(i+j+1)]_s}\right) \\
& = &
q^{-n/4} s^{-n^2/2} (-1)^n
\left(\prod_{0 \le j < i < n} {[3(i-j)]_s \over 3[i-j]_s}\right)
\left(\prod_{i=0}^{n-1} {\prod_{j=1}^{3i+1} [j]_s \over
\prod_{j=1}^{n+i} [j]_s}\right)
\end{eqnarray*}
Note that the second factor is Andrews' $q$-enumeration of descending plane
partitions~\cite{Robbins:intelligencer} with $q$ replaced by $s$. Taking the
limit as $\epsilon \to 0$ and combining with equation~(\ref{ean}), the
factors of $q$ and $-1$ cancel, the factors of $s$ become factors of $1$, and
the brackets disappear.  The result is
$$1!4!7!\ldots(3n-2)! \over n!(n+1)!(n+2)!\ldots(2n-1)!.$$
This completes the proof of Theorem~\ref{thmain}.

For general $x$, the matrix $M$ becomes
$$M_{i,j} = {x^2-4x \over s^{i+j+1} + 2 - x + s^{-(i+j+1)}}.$$
There are two other values of $x$ when the denominator is a cyclotomic
(Laurent) polynomial in some power of $s$, namely $x=2$ and $x=3$. In the
former case, $-M/4 = S(n;s^2,s^4)$, whose determinant is given by
Lemma~\ref{lsdet}; alternatively, the determinant may also be derived from
Lemma~\ref{ltdet}. In the latter case, $-M/3 = S'(n;s,s^3)$, where
$$S'(n;s,t)_{i,j} = {s^{i+j+1} + 1 \over t^{i+j+1} + 1}.$$
A variation of Lemma~\ref{lsdet} establishes the determinant of $S'(n;s,t)$;
the leading coefficient is simply the determinant of $S(n;t,t^2)$. These
manipulations clearly lead to product formulas for $A(n;2)$ and $A(n;3)$, and
in particular, to a proof of Theorem~\ref{thmain3}.  We omit
the details of rearranging and cancelling factors to put the product
formulas in their standard form.

Finally, we use Theorem~\ref{thkorepin} to prove Theorem~\ref{thx}.
Recall the variables $x_i$ and $y_i$
in the definition of $Z$, which are not to be confused with the $x$
of $A(n;x)$.  If we set $x_i = \frac12 + f_i\epsilon$
and $y_j = f_j\epsilon$ for some $f_i$'s such that $f_{n-1-i} = -f_i$,
then $Z(n;X,Y)$ again converges to $A(n;x)$ up to normalization
as $\epsilon \to 0$.  In this case,
the corresponding matrix $M$ is given by
$$M_{i,j} = {x^2 - 4x \over s^{f_i-f_j} + 2 - x + s^{f_j-f_i}}.$$
This matrix $M$ possesses the symmetry $(i,j) \mapsto (n-1-i,n-1-j)$, \ie, it
commutes with the antidiagonal permutation matrix $P$. Therefore, a change of
basis divides $M$ into blocks corresponding to the eigenspaces of $P$. 
Therefore the determinant of $M$ is the product of the determinants of the
blocks. This is the origin of the factorization of the $A(n;x)$'s
into the $B(n;x)$'s.

%\nocite{Zeilberger:asm2}
%\bibliographystyle{plain}
%\bibliography{pp}

\end{document}